\documentclass{amsart}
\usepackage{amssymb}
\usepackage{amsfonts}
\usepackage{amsmath}
\usepackage[legalpaper,bookmarks=true,colorlinks=true,linkcolor=blue,citecolor=blue]{hyperref}
\usepackage{graphicx}%
\usepackage{fancyhdr}
\usepackage{color}
\usepackage[mathlines]{lineno}
\usepackage{lscape}
\usepackage{epsfig}
\usepackage{natbib}
\usepackage{geometry}
\usepackage{tgbonum}
\usepackage[]{listings}

\fontfamily{qcr}\selectfont

\newtheorem{theorem}{Theorem}
\theoremstyle{plain}

\newtheorem{definition}{Definition}

\newtheorem{lemma}{Lemma}

\newtheorem{fact}{Fact}

\numberwithin{equation}{section}

\newcommand{\Bin}{\bigskip \noindent}

\newcommand{\Ni}{\noindent}

\begin{document}
\Large
\title{Asymptotics of summands I: square integrable independent random variables}
\author{Aladji Babacar Niang}
\author{Gane Samb Lo}
\author{Moumouni Diallo}

\begin{abstract} This paper is part of series on self-contained papers in which a large part, if not the full extent, of the asymptotic limit theory of summands of
independent random variables is exposed. Each paper of the series may be taken as review exposition but specially as a complete exposition expect a few exterior resources. For graduate students and for researc\emph{h}ers (beginners or advanced), any paper of the series should be considered as a basis for constructing new results. The contents are taken from advanced books but the organization and the proofs use more recent tools, are given in more details and do not systematically follow previous one. Sometimes, theorems are completed and innovated.\\

\noindent $^{\dag}$ Aladji Babacar Niang\\
LERSTAD, Gaston Berger University, Saint-Louis, S\'en\'egal.\\
Email: niang.aladji-babacar@ugb.edu.sn, aladjibacar93@gmail.com\\

\noindent $^{\dag \dag}$ Gane Samb Lo.\\
LERSTAD, Gaston Berger University, Saint-Louis, S\'en\'egal (main affiliation).\newline
LSTA, Pierre and Marie Curie University, Paris VI, France.\newline
AUST - African University of Sciences and Technology, Abuja, Nigeria\\
gane-samb.lo@edu.ugb.sn, gslo@aust.edu.ng, ganesamblo@ganesamblo.net\\
Permanent address : 1178 Evanston Dr NW T3P 0J9,Calgary, Alberta, Canada.\\

\noindent $^{\dag \dag \dag}$ Dr Moumouni Diallo\\ 
Université des Sciences Sociale et de Gestion de Bamako ( USSGB)\\
Faculté des Sciences Économiques et de Gestion (FSEG)\\
Email: moudiallo1@gmail.com\\




\end{abstract}
\maketitle
\newpage
\tableofcontents
\newpage
\section{Introduction}

\Ni The largest part of the asymptotic theory of partial sums of random variables concentrated on independent random variables over at least two centuries. Almost all the greatest scientists in probability Theory (L\'evy, Kolmogorov, Lyapounov, Lynderberg, Gnedenko, Feller, etc.) engaged themselves in such an enterprise. Besides, a very large part of the current theory on dependent sequence of random variables is based on transformations of independence structures, for example on notions of \textit{nearness} of the dependence to independence (weak dependence, $\phi$-mixing, associated sequence, independent increments, etc.).\\

\Ni So it is important to have the deepest knowledge of that past. In \cite{ips-mfpt-ang}, we introduced some important elements of that theory (Central limit theorems, laws of the large numbers, law of the iterated logarithm, zero-one laws, etc.).\\

\Ni We are beginning a series on self-contained papers in which a large part of the \textit{central limit theorem}, if not the full extent, of the asymptotics of summands of independent random varianles will exposed. Each paper of the series may be taken as review exposition but specially as a complete exposition expect a few exteroir ressources. For graduate students and for researchers (beginners or advanced), any paper of the series should be considered as a basis for constructed new results. The contents are taken from advanced books but the organization and the proofs used more recent tools, are given in more details and do not systematically follow previous one. Sometimes, theorem are completed and innovated.\\

\Ni In this first paper of the series, we focus of the full characterization of the \textit{CLT} problem of independent summands for square integrable random variables. The main material is extracted from \cite{loeve} as a general guide. But we use arguments from our previous works (\cite{ips-mestuto-ang}, \cite{ips-probelem-ang}, \cite{ips-mfpt-ang}, \cite{ips-wcia-ang}, etc.) to have unified and a self-contained \cite{ips-mfpt-ang}. In particular, the text on the weak convergence of bounded measures and its exppression on $\mathbb{R}^k$ ($k\geq 1$) povides tools to make the conclusions in \cite{loeve} more clear, as we expect as least.\\

\Ni Papers of this series focus on complete mathematical texts rather than on a historical review of contributions of differents authors. We refer to \cite{loeve} for that aspect.\\  

\Ni Since the set of weak limits of independent summands for square integrable coincide with the set of infinitely decomposable laws, this paper will deal with the most important properties of such probability laws to the extent needed here. More developments, including the study of stable laws, will be given in the other papers of the series.\\

\Ni Let us introduce the problem, after we provide some notations. below, all sequences and all arrays of random variables have all their elements defined on a same probability space. So, we set a generic probability space
 $(\Omega, \mathcal{A}, \mathbb{P})$.\\

\noindent Following \cite{loeve}, we transform the study of sums of independent and centered random variables 

$$
S_n=X_1+...+X_n, \ n\geq 1, 
$$

\Bin (with the notation $\sigma_i^2=\mathbb{E}X_i^2$, $i\geq 1$, if they exist), by summands

$$
S_n=\sum_{1\leq k \leq k_n} X_{k,n}, \ n\geq 1,
$$

\Bin where for each $n\geq 1$, the family $\{X_{k,n}, \ 1 \leq k \leq k_n=k(n)\}$ is a family of independent and centered random variables such that $F_{k,n}$ stands for the \textit{cdf} of $X_{k,n}$ and $\sigma_{k,n}^2=\mathbb{E}X_{k,n}^2$, $1\leq k \leq k_n$. We suppose also that $k_n \rightarrow +\infty$ as $n\rightarrow +\infty$.\\

\Ni \textbf{Notations}. The notation already given and completed by the notation $f_{kn}$ for the characteristic function of $X_{k,n}$, are fixed for once.\\ 

\noindent In the case of simple summands, we have for each $n\geq 1$, $k_n=n$ and $X_{k,n}=X_k$ for $k \in [1,n]$. Here, the rows $(X_{k,n})_{1\leq k \leq k_n}$ are such that each of them is obtained by adding one element to the predecessor. But, in the general case, no relation between families 
$\mathcal{E}_n=\{X_{k,n}, \ 1 \leq k \leq k(n)\}$ is required. Also, in the case of the simple sequence $(X_k)_{k\geq 1}$, the studied array for each $n\geq 1$, is
 $\{X_1/s_n, \cdots, X_n/s_n\}$ where $s_n^2=\mathbb{V}ar(X_1+\cdots+X_n)$.

\noindent Here, we are going to investigate the general problem of finding all the possible weak limits do $S_n$. Without restrictions, we this may lead to trivial results. So we have to fix a general frame in which the study will be done. In doing, the best way seems to go back to the complete theory of Feller-Levy-Lynderberg and there, discover the following two fundamental hypotheses.\\

\section{The Bounded Variance Hypothesis (\textit{BVH}) and the Uniformly Asymptotic Negligibility (\textit{UAN})} \label{03_id_sec_01}

\Ni  Given a sequence $(X_k)_{k\geq 1}$ of independent, centered and square-integrable random variables, we set $\sigma_k^2=\mathbb{E}X_k^2$, $k\geq 1$, for $n\geq 1$

$$
s^2_n=\sum_{k=1}^{n} \sigma_k^2,  \ \ t_n^2=\max \{\sigma_k^2, \ \ 1\leq k \leq n\} \ \ and \ \ B_n=t_n^2/s_n^2,
$$

\Bin $k(n)=n$ for each $n\geq 1$ and for each $n\geq 1$

$$
\{X_{k,n}, \ 1\leq k\leq k(n)\}=:\left\{\frac{X_k}{s_n}, \ 1\leq k\leq k(n)\right\}
$$

\Bin and

$$
S_n=\sum_{k=1}^{k(n)} X_{k,n}.
$$

\Bin The Feller-Levy-Lynderberg (\textit{F2L}) theorem (see \cite{loeve}, or \cite{ips-mfpt-ang}, Chapter 7, Section 2, Part B) ensures that: \label{HVB-UAN}\\

\begin{equation}
S_n \rightsquigarrow \mathcal{N}(0,1) \ \ [WC] \ \ and \ \  B_n\rightarrow 0 \ \ [NG] \label{Ng}
\end{equation}

\Bin if and only if, for any $\varepsilon>0$,

\begin{equation}
L_n(\varepsilon)=\sum_{1\leq k \leq k(n)} \int_{|X_{k,n}|\geq \varepsilon} X_{k,n}^2 \ d\mathbb{P} \rightarrow 0. \label{lyndeberg}
\end{equation}

\Bin Let us see how behave the two following important quantities in that frame :

\begin{equation}
U(n,\varepsilon)= \sup_{1\leq k \leq n} \mathbb{P}(|X_{k,n}|\geq \varepsilon) \ \ and \ \  MV(n)=\sum_{1\leq k \leq k(n)} \mathbb{V}ar(X_{k,n}).
\end{equation}

\Bin We have, by Markov inequality,

\begin{equation}
U(n,\varepsilon)\leq \sup_{1\leq k \leq n} \frac{\mathbb{E} X_{k,n}^2}{\varepsilon^2 s_n^2}=\varepsilon^{-2} B_n, \label{UAN2}
\end{equation}

\noindent and $MV(n)=1$ for all $n\geq 1$ and for $c=1$

$$
\sup_{n\geq 1} MV(n)=c<+\infty.
$$

\Bin The theory we are going to develop in a more general case needs the properties we just introduced with specific names.

\Bin \textbf{Definition}. Under the notation given above, we say that :\\

\Ni (i) the \textbf{Bounded Variance Hypothesis (\textit{BVH})} holds if and only if

$$
\sup_{n\geq 1} MV(n)=c<+\infty;
$$

\Bin (ii) the \textbf{Variance Convergence Hypothesis (\textit{VCH})} holds if and only if

$$
MV(n) \rightarrow c \in ]0,+\infty[;
$$

\Bin (iii) the \textbf{Uniformly Asymptotic Negligibility (\textit{UAN})} holds if and only if, for any $\varepsilon>0$,

$$
U(n,\varepsilon) \rightarrow 0.
$$

\Bin We express the \textit{F2L} theorem as follows: If \\

\noindent (a) the \textit{(UAN)} condition holds;\\

\noindent (b) the \textit{(BVH)} holds, then 

$$
S_n \rightsquigarrow \mathcal{N}(0,1) \ \ \text{\textit{if \ and \ only \ if}} \ \ L_n(\varepsilon)\rightarrow 0, \ \text{\textit{for  any}} \ \varepsilon>0.
$$

\Bin In this particular case, the characteristic function of the weak law $\mathcal{N}(0,1)$

$$
\psi_{\infty}(u)= \exp(-u^2/2)
$$

\Bin is such that, for any $p \in \mathbb{N}$, $p>0$, $\psi_{\infty}^{1/p}$ defined by

$$
\psi_{\infty}^{1/p}(u)= \exp(-p^{-1}u^2/2)
$$

\Bin is still a characteristic function, actually of a $\mathcal{N}(0, \ p^{-1})$ law. Let us denote by $\mathcal{C}_f$ the class of all characteristic functions
$f : \mathbb{R}\rightarrow \mathbb{C}$ on $\mathbb{R}$.\\

\Bin Now we may set our general task to be done: (\textit{Task}). Given the \textit{UAN} and the \textit{BVH} conditions, what is the class of all possible limits 
$Z_{\infty}$ of characteristic function $\psi_{\infty}$. By the particular case of \textit{F2L} theorem, we may think that the searched class can be

$$
\mathcal{C}_{fid}=\{\psi \in \mathcal{C}_f: \ \ (\forall p>1), \ \psi^{1/p} \in \mathcal{C}_f\}.
$$

\Bin We define $\mathcal{C}_{fid}$ as the class of infinitely divisible characteristic functions. At least, in the current version of the central limit theorem, the Gaussian weak limit is in $\mathcal{C}_{fid}$.\\

\Ni We are going to see that the suggestion in the description of the \textit{task} is effectively the global solution.\\

\Ni In the sequel, we will devote Section \ref{03_id_sec_02} on infinitely divisible (or decomposable) laws. In Section \ref{03_stable_GEN_03}, we finish the 
\textit{task} we have given to ourselves under the \textit{UAN} Condition and the \textit{BVC}. Finally in \ref{03_Gauss-Poisson}, rediscover the characterization of the \textit{CLT} to a Gaussian law and that of the \textit{CLT} to a Poisson law. In the next element of the series, we proceed to a general theory with non-necessarily square-integrable random variables.

\section{Class of infinitely divisible (or decomposable) laws on $\mathbb{R}$} \label{03_id_sec_02}

\subsection{Definitions and examples} $ $\\

\Ni The basic definition is the following.

\begin{definition} \label{idf_def} A characteristic function $\psi \in \mathcal{C}_{f}$ is infinitely decomposable (\textit{idecomp}), denoted by $\psi \in \mathcal{C}_{fid}$ if and only if for all positive integer $p$, $\psi^{1/p}$ is still a characteristic function.
\end{definition}
 
\Bin Let us explain the notion of \textit{idecomp} in terms of random variables. Suppose that $\psi \in \mathcal{C}_{f}$, $p\geq 1$ and $\psi_p=\psi^{1/p}$. Suppose that $\psi_p$ is the characteristic function (\textit{ch.f}) of a probability measure $\mathbb{P}_{p}$. By the Kolmogorov theorem, it is possible to construct a probability space $(\Omega_p, \mathcal{A}_p,\mathbb{P}_p)$ holding independent real-valued random variables $Z_{p}$, $Z_{1,p}$, $\cdots$, $Z_{p,p}$ having all the \textit{cha.f} $\psi_p$, i.e,

$$
\psi_p(u)= \mathbb{E}_{\mathbb{P}_{p}} \exp(iuZ_{p})  = \int \exp(iux) \ d\mathbb{P}_{p}.
$$

\Bin It is clear that $\psi_p^p$ is the \textit{cha.f} of $S_p=Z_{1,p}+\cdots+Z_{p,p}$. As well, $\psi$ is the \textit{cha.f} of a probability measure $\mathbb{P}_{\psi}$ on $\mathbb{R}$ and let us denote by $Z_{\psi}$ a random variable with $\mathbb{P}_{\psi}$ as probability law.\\

\Bin We easily see that we may extend the definition as follows. In the definition below, we use the notion of \textit{idecomp} probability law at the place of \textit{idecomp characteristic function} or \textit{idecom random variable}.

\begin{definition} \label{idf_defExt} Let $Z$ be a real-valued random variable with probability law $\mathbb{P}_Z$ and \textit{ch.f} $\psi_Z$. $\mathbb{P}_Z$ is \textit{idecom} (equivalently $\psi_Z$ is idecomp or $Z$ is idecomp) if and only if one of the following assertions holds:\\

\noindent (i) For all $p \in \mathbb{N}\setminus \{0\}$, $\psi_Z^{1/p}$ is a \textit{cha.f}.\\

\noindent (ii) For all $p \in \mathbb{N}\setminus \{0\}$, there exists a \textit{cha.f} $\psi_p$ such that $\psi_Z=\psi_{p}^{p}$.\\

\noindent (iii) For all $p \in \mathbb{N}\setminus \{0\}$, there exists a probability $\mathbb{P}_{p}$ on $\mathbb{R}$ such that

$$
\mathbb{P}_Z=\mathbb{P}^{\otimes p}_{p},
$$

\Bin that is, $\mathbb{P}_{Z}$ is the convolution product of $\mathbb{P}_{p}$ by itself $p$ times.\\

\Ni (iv) For all $p \in \mathbb{N}\setminus \{0\}$, there exists a sequence $Z_{1,p}$, $\cdots$, $Z_{p,p}$ of independent and identically distributed real-valued random variables such that

$$
Z =_{d} Z_{1,p}+\cdots+Z_{p,p}.
$$
\end{definition} 

\Bin \textbf{Examples}. Let us give some quick examples.\\

\Ni \textbf{Example 1}. (Degenerate random variable). Let $Z=a$, \textit{p.s} of \textit{cha.f}

$$
\psi_Z(t)=e^{iat}, \ t \in \mathbb{R}.
$$

\Bin For $p\geq 1$, $\psi_Z(t)^{1/p}=e^{i(a/p)t}$, which is the \textit{cha.f} of the degenerate \textit{r.v} $Z_{p}=a/p$.\\

\Ni \textbf{Example 2}. (Gaussian random variables). Let $Z \sim \mathcal{N}(m,\sigma^2)$, $m \in \mathbb{R}$, $\sigma \in \mathbb{R}_{+}\setminus \{0\}$ with \textit{cha.f}

$$
\psi_Z(t)=\exp (imt - \sigma^{2}t^2/2), \ t \in \mathbb{R}.
$$

\Bin For $p\geq 1$, we have 

$$
\psi_Z(t)^{1/p}=\exp (i(m/p)t - (\sigma/\sqrt{p})^2t^2/), \ t \in \mathbb{R},
$$

\Bin which is the \textit{cha.f} of a $\mathcal{N}(m/p,\sigma^2/p)$ \textit{r.v}.\\

\Ni \textbf{Example 3}. (Translated Poisson random variables). Let $Z \sim \mathcal{P}(a,\lambda)\equiv + \mathcal{P}(\lambda)$, $a \in \mathbb{R}$, $\lambda \in \mathbb{R}_{+}\setminus \{0\}$ with \textit{cha.f}

$$
\psi_Z(t)=\exp \left(iat + \lambda(e^{it}-1) \right).
$$

\Bin For $p\geq 1$, we have 

$$
\psi_Z(t)^{1/p}=\exp \left(i(a/p)t + (\lambda/p)(e^{it}-1) \right), \ t \in \mathbb{R},
$$

\Bin which is the \textit{cha.f} of a $\mathcal{P}(a/p,\lambda/p)$ \textit{r.v}.\\

\Ni \textbf{Example 4}. (Gamma random variables). Let $Z \sim \gamma(a,b)$, $a>0$, $b>0$, with \textit{cha.f}

$$
\psi_Z(t)=\left(1 - it/b\right)^{-a}, \ t \in \mathbb{R}.
$$

\Bin For $p\geq 1$, we have 

$$
\psi_Z(t)^{1/p}=\left(1 - it/b\right)^{-(a/p)}, \ t \in \mathbb{R},
$$

\Bin which is the \textit{cha.f} of a $\gamma(a/p, b)$ \textit{r.v}.\\

\Ni \textbf{Example 5}. (Cauchy random variables). Let $Z \sim Ca(a,b)$, $a\in \mathbb{R}$, $b>0$ with \textit{cha.f}

$$
\psi_Z(t)=\exp(iua -b|t|), \ t \in \mathbb{R}.
$$

\Bin For $p\geq 1$, we have 

$$
\psi_Z(t)^{1/p}=\exp(iu(a/p) -(b/p)|t|), \ t \in \mathbb{R},
$$

\Bin which is the \textit{cha.f} of a $Ca(a/p, b/p)$ \textit{r.v}.\\

\Bin Now, let us focus on properties of such laws.\\

\subsection{Properties} $ $\\

\Ni \textbf{Property 1}. \label{property_01} If $\psi_1$ and $\psi_2$ are two \textit{idecomp} \textit{cha.f}, then $\psi=\psi_1 \psi_2$ is an \textit{idecomp} \textit{cha.f}.\\

\Ni \textbf{Proof}. Suppose that $\psi_1$ and $\psi_2$ are two \textit{idecomp} \textit{cha.f} and let $p\geq 1$. Thus $\psi^{1/p}=\psi_1^{1/p} \psi_2^{1/p}$ is the \textit{cha.f} of the convolution product of the probability measures associated to the \textit{cha.f} $\psi_i^{1/p}$ ($i \in \{1,2\}$).\\

\Ni \textbf{Property 2}. \label{property_02} If $\psi$ is an \textit{idecomp} \textit{cha.f}, the conjugate $\overline{\psi}$ is also an \textit{idecomp} \textit{cha.f} and the complex square norm $\|\psi\|^2$ is an \textit{idecomp} \textit{cha.f}.\\

\Ni \textbf{Proof}. Let $\psi$ be the \textit{cha.f} of $X$, i.e., $\psi(t)=\psi_{X}(t)=\mathbb{E}(e^{itX})$, it is clear that

$$
\mathbb{E}(e^{-itX})=\mathbb{E}(\overline{e^{itX}})=\overline{\mathbb{E}(e^{itX})}=\overline{\psi}(t).
$$

\Bin This and $\mathbb{E}(e^{-itX})=\psi_{-X}(t)$ for $t\in \mathbb{R}$ show that $\psi$ is a \textit{cha.f}. It is also direct to see that $X$ and $-X$ are \textbf{idecomp} or {non-idecomp} at the same time or not. Finally by Property 1, $\|\psi\|^2=\psi \overline{\psi}$ is \textit{idecomp} if $\psi$ is.\\

\Ni \textbf{Property 3}. \label{property_03} If $\psi$ is an \textbf{idecomp} \textbf{cha.f}, then $\psi^{1/n}$ converges to 1 everywhere, as \ $n\rightarrow +\infty$.\\

\Ni \textbf{Proof}. Suppose that $\psi$ is an \textit{idecomp} \textit{cha.f}. Let us denote, for all $n\geq 1$, $\psi_n=\psi^{1/n}$, that is a \textit{cha.f}. But
$\|\psi\|\leq 1$ and $\|\psi_n\|^2=\|\psi\|^{2/n}$ converges to $g$ with $g=0$ on $\psi=0$ and $g=1$ on $\psi\neq 0$. Let us show that  $\psi$ cannot take the null value. Indeed $\psi$ is continuous (at zero in particular) and $\psi(0)=1$. So $\psi>1/2$ on an interval $]-r,r[$, $r>0$ and next $g=1$ on $]-r,r[$. But the function
$h\equiv 1$ is the \textit{cha.f} of the random variable $Z=0$. By Proposition in \cite{billingsley} (see page 388), we get that $g=h$ and then $g=1$ everywhere, so $\|\psi_n\|^2 \rightarrow 1$. This ensures that $\psi$ does not take the null value. Finally, we get rid of the norm by

$$
\psi^{1/n}=\exp\left(\frac{1}{n} \log \psi \right) \rightarrow 1 \ as \ n\rightarrow +\infty. 
$$

\Ni \textbf{Property 4}. \label{property_04} Let $(\psi_n)_{n\geq 1}$ be a sequence of \textit{idecomp} \textit{cha.f}'s such that $\psi_n \rightarrow \psi$ and $\psi$ is continuous at zero. Then $\psi$ is an \textit{idecomp cha.f}.\\ 

\Ni \textbf{Proof}. Let $\mathcal{C}_{fid} \ni \psi_p \rightarrow \psi$ and $\psi$ is continuous at zero. For any fixed $q\geq 1$, $|\psi_p|^{2/q} \rightarrow |\psi|^{2/q}$. Since the $|\psi_p|^{2/q}$, are \textit{cha.f} and $|\psi^{2/q}|$ is continuous at zero, it comes that $\psi|^{2/q}$ is a \textit{cha.f} for any $q\geq 1$. So $|\psi|^2$ is an \textit{idecomp} \textit{cha.f} and by property 3, $\psi$ is nonwhere zero and next

$$
\psi_p^{1/q}=\exp\left(\frac{1}{q} \log \psi_p \right) \rightarrow \exp\left(\frac{1}{q} \log \psi \right)=\psi^{1/q} 
$$

\Bin is a \textit{cha.f} by the Levy continuity theorem.\\

\Ni \textbf{Property 5}. \label{property_05} An  \textit{cha.f} $\psi$ is \textit{idecomp} if and only if it is limit of a sequence of products Poisson type \textit{cha.f}.\\

\Ni \textbf{Proof}. If $\psi$ is a limit of a sequence of products Poisson type \textit{cha.f}, it is \textit{idecomp} by Property 4, since products of Poisson type \textit{cha.f} are \textit{idecomp} \textit{cha.f}.\\

\Ni Conversely, let us be given an \textit{idecomp} \textit{cha.f} $\psi$. Since $\psi$ is non-where equal to zero (Property 3), we have

$$
\log \psi = \lim_{p\rightarrow +\infty} p(\psi^{1/p}-1).
$$

\Bin For $p\geq 1$, let us denote by $F_{p}$ the \textit{cdf} associated with the \textit{cha.f} $\psi_p=\psi^{1/p}$. So we have

$$
\Psi_p(t)=p(f^{1/p}-1) =\int p\left(e^{itx}-1\right) \ dF_p(x).
$$

\Bin Since the function $\hookrightarrow p\left(e^{itx}-1\right)$ is bounded on $\mathbb{R}$, it is locally integrable and $\lambda_{F_p}$ is a finite measure, we may apply Lebesgue Dominated theorem and we can conclude that for any fixed $p\geq 1$,

$$
\Psi_p(t)=p(f^{1/p}-1) =\lim_{0<a\rightarrow +\infty} \int_{-a}^{a} p\left(e^{itx}-1\right) \ dF_p(x)=:\lim_{0<a\rightarrow +\infty} \Psi_{p,a}.
$$

\Bin By continuity of the integrand, the integral $\Psi_{p,a}$ is limit of Riemann-Stieltjes, which are of the form

$$
\sum_{1\leq j\leq k(p,a)} pb_{j,p} \left( e^{ic_{j,p}u}-1 \right),
$$ 

\noindent which are sums of logarithms of Poisson type \textit{cha.f}. Hence $\exp(\Psi_{p,a})$ are \textit{cha.f} and next $\exp(\Psi_{p})$ is a \textit{cha.f} as limit of the sequence $\exp(\Psi_{p,a})$.\\

\Ni Finally $\psi$ is limit of \textit{cha.f} of the form $\exp(\Psi_{p})$, which is  a sequence of products of Poisson type \textit{cha.f}.\\

\Ni \textbf{Property 6}. \label{pageprop4} A \textit{cha.f} is \textit{idecomp} if and only if it is limit of a sequence of products of \textit{cha.f} of Poisson type laws.\\

\Ni \textbf{Proof} Let $\psi$ be a \textit{cha.f}. Let $\psi_p$ a product of \textit{cha.f} of type Poisson laws

$$
\psi_p(t)=\prod_{j=1}^{k(p)} \exp\left(ia_{j,p}t + b_{j,p} \left(e^{ic_{j,p}t}-1\right) \right), \ t \in \mathbb{R},
$$

\Bin where the $a_{j,p}$'s and $c_{j,p}$'s are real numbers and the $b_{j,p}$'s positive numbers. We have, for $p$ fixed and for $q\geq 1$

$$
\psi_p(t)^{1/q}=\exp\left(i \left\{\frac{1}{q}\sum_{j=1}^{k(p)} a_{j,p}\right\} t + \frac{1}{q} \sum_{j=1}^{k(p)} b_{j,p} \left( e^{i c_{j,p} t}-1\right) \right), \ \ t \in \mathbb{R}.
$$

\Bin This is still a product of \textit{cha.f}'s of type Poisson type laws and is a \textit{cha.f}. If $\psi_p \rightarrow \psi$, thus by Property 4, $\psi$ is \textit{idecomp}.\\

\Ni We will need more facts on \textbf{cha.f}'s that we will introduce when needed,\\

\Ni We begin by studying the case of bounded variances. First, we deal with three important results that constitute the pillars of the current theory. \\

\subsection{The three pillars of that theory} \label{03_03_ssec_lem_com}

\Ni In this subsection, we assume that both the \textit{UAN} and the \textit{BVH} hold. 

\begin{lemma} \label{03_03_lem_01} (Comparison Lemma) The complex function $\log f_{k,n}$ is well-defined and for any $u\in \mathbb{R}$

$$
\sum_{k=1}^{k(n)}\left\{\log f_{k,n}(u) -(f_{k,n}(u)-1)\right\} \rightarrow 0,
$$

\Ni as $n\rightarrow +\infty$.
\end{lemma}

\Ni \textbf{Proof}. Let $u\in \mathbb{R}$ fixed and $n\geq 1$. Then for any $k\in \{1,\cdots,k(n)\}$, we have the one order expansion

\begin{equation*}
f_{k,n}(u)=1+ \theta_{k,n} u^2 \sigma_{k,n}^2/2,
\end{equation*}

\Ni with $|\theta_{k,n}|<1$ and $|\circ|$ stands for the norm in $\mathbb{C}$ or the absolute value when applied to real numbers. In all this chapter, numbers of the form $\theta_{\circ}$, possibly written with primes or double primes, are only required to have norms less than one and their values are not important. So, we get

\begin{equation*}
\max_{1\leq k \leq k(n)} |f_{k,n}(u)-1|\leq \frac{u^2 B_n}{2} \rightarrow 0 \ as \ \ n\rightarrow +\infty.
\end{equation*}

\Ni Next for $v_{k,n}=\theta_{k,n} u^2 \sigma_{k,n}^2/2$, we surely have that $\max_{1\leq k \leq k(n)} |v_{k,n}| \leq (u^2 B_n)/2$ goes to zero. We also have for all $u\mathbb{R}$,

\begin{eqnarray*}
\log f_{k,n}(u)&=& \log (1 + (f_{k,n}(u)-1))=\log (1 + v_{k,n})=v_{k,n} + \theta_{k,n}^{\prime} v_{k,n}^2\\
&=& (f_{k,n}(u)-1) + \theta_{k,n}^{\prime} v_{k,n}^2,
\end{eqnarray*}

\Ni which leads to, as $n\rightarrow +\infty$,

\begin{eqnarray}
\left| \sum_{k=1}^{k(n)}\left\{\log f_{k,n}(u) -(f_{k,n}(u)-1)\right\} \right| &\leq& \sum_{1}^{k(n)} |\theta_{k,n}^{\prime}| v_{k,n}^2 \label{boundX_01}\\
&\leq& \sum_{k=1}^{k(n)} \frac{u^4}{4} |\theta_{k,n}|^2 \sigma_{k,n}^4 \ \ \ \ (L3)\notag\\
&\leq& \frac{u^4 B_n}{4} \sum_{k=1}^{k(n)} \sigma_{k,n}^2 \ \ \ \ (L4)\notag \\
&\leq& \frac{c u^4 B_n}{4} \rightarrow 0. \notag
\end{eqnarray}

\Bin The proof of Lemma \ref{03_03_lem_01} is over. $\square$\\

\Ni Now, let us use new expressions of the results in Lemma \ref{03_03_lem_01}. Since the variables $X_{k,n}$ are centered, we have

\begin{equation*}
\forall n\geq 1, \  \forall 1\leq k\leq k(n), \ \int X_{k,n} \ d\mathbb{P}=\int x \ dF_{k,n}(x)=0 \ and \ \ \int x \ \ dF_{k,n}(x)=\sigma^{2}_{k,n}. 
\end{equation*}

\Bin Let us set, for  $n\geq 1$,

\begin{equation*}
\psi_n(u)\equiv \sum_{k=1}^{k(n)}(f_{k,n}(u)-1)=\sum_{k=1}^{k(n)} \int \left(e^{iux}-1\right) \ dF_{k,n}(x), \ u\in \mathbb{R}.
\end{equation*}

\Bin By using the remark that $\mathbb{E}X_{k,n}=0$, i.e. $\int x \ dF_{k,n}(x)=0$, we get

\begin{eqnarray*}
\psi_n(u)&=&\sum_{k=1}^{k(n)} \int \left(e^{iux}-1-iux\right) \ dF_{k,n}(x)\\
&=& \int \left(e^{iux}-1-iux\right) \ \sum_{k=1}^{k(n)} dF_{k,n}(x)\\
&=& \int \frac{1}{x^2} \left(e^{iux}-1-iux\right) \  x^2 \sum_{k=1}^{k(n)} dF_{k,n}(x), \ u \in \mathbb{R}.
\end{eqnarray*}

\Ni But, by putting

$$
dK_n(x)=x^2 \sum_{k=1}^{k(n)} dF_{k,n}(x),
$$

\Bin we get

\begin{equation*}
\psi_n(u)= \int \frac{1}{x^2} \left(e^{iux}-1-iux)\right) \ dK_n(x).
\end{equation*}

\Bin Finally,  Lemma \ref{03_03_lem_01} can be expressed as \\

\begin{lemma} \label{03_03_lem_02}
\begin{equation*}
\forall u \in \mathbb{R}, \log\left(\prod_{k=1}^{k(n)} f_{k,n}(u)\right) - \psi_n(u)\rightarrow 0, \ as \  n\rightarrow +\infty,
\end{equation*}

\Bin where

$$
dK_n(x)=x^2 \sum_{k=1}^{k(n)} dF_{k,n}(x)
$$

\Ni and

\begin{equation*}
\psi_n(u)= \int \frac{1}{x^2} \left(e^{iux}-1-iux\right) \ dK_n(x).
\end{equation*}
\end{lemma}

\Bin This lemma becomes the second pillar. The third is the following

\begin{lemma} \label{03_03_lem_03}
For any $n\geq 1$, $\exp(\Psi_n)$ is an \textit{idecomp} \textit{cha.f} and is the \textit{cha.f} of a centered random variable of variance

$$
\int \ dK_n(x)=s_n^2.
$$
\end{lemma}

\Bin \textbf{Proof}. Let $n\geq 1$ be fixed. We have

$$
\Psi_n(u)=\int \ g(u,x) \ dK_n(x) \ with \ g(u,x)=\frac{e^{iux}-1-iux}{x^2}, \ x \in \mathbb{R}.
$$

\Bin Clearly $g$ is continuous on $\mathbb{R} \times \mathbb{R}^{\ast}$ (with $\mathbb{R}^{\ast}=\mathbb{R}\setminus \{0\}$) and for $u$ fixed, $g(u,0)$ is the extension of $g(u,x)$ by limit, since an expansion at zero gives

$$
g(u,x)=\frac{1+iux-u^2x^2/2-1-iux+O(x^3)}{x^2} \rightarrow -u^2/2 \ as \ x\rightarrow 0.
$$

\Bin So, for $u$ fixed, $x\mapsto g(u,x)$ is continuous everywhere. Moreover we have

\begin{equation}
\forall u\in \mathbb{R}, \ \forall x\in \mathbb{R}^{\ast}, \ |g(u,x)|\leq \frac{2}{x^2} + \frac{|u|}{|x|} \label{bound_01}
\end{equation}

\Bin and

$$
\int x^{-2} \ dK_n(x)=\sum_{k=1}^{k(n)} \int dF_{k,n}(x)=k(n),
$$

\Bin and, by using $(|x|\leq 1 +x^2)$

\begin{eqnarray*}
\int |x|^{-1} \ dK_n(x)&=&\sum_{k=1}^{k(n)} \int \frac{x^2}{|x|} \ dF_{k,n}(x)\\
&=&\sum_{k=1}^{k(n)} \int |x| \ dF_{k,n}(x)\\
&\leq&\sum_{k=1}^{k(n)} \int (1+x^2)\ dF_{k,n}(x)\\
=k(n)+s_n^2. 
\end{eqnarray*}

\Bin We conclude that $g(u,x)$ is bounded by $g_0(x)=2x^{-2}+|ux^{-1}|$ which is $K_n$-integrable. So by the dominated convergence theorem, $\Psi_n$ is continuous at zero. Also, as an improper Riemann-Stieltjes integral, for $\varepsilon>0$ fixed, we can find a number $A>0$ such that for $a\geq A$

$$
\left| \psi_n(u) - \Psi_{n,a}(x)\right|<\varepsilon \ with \ \Psi_{n,a}(u)=\int_{-a}^{a} g(u,x) \ dK_n(x).
$$

\Bin Now, since $\Psi_{n,a}(u)$ is continuous, it a limit of a sequence of  of Riemann-Stieltjes sums: there exists a partition of $[-a,a]$

$$
-a=x_{0,p}<\cdots <x_{j-1,p}<x_{j,p}<\cdots<x_{\ell(p),p}=a
$$

\Bin and a sequence of points $c_{j,p} \in (x_{j,p}, x_{j+1,p})$, $0\leq j\leq \ell(p)-1$, 

$$
S_p(u)=\sum_{j=0}^{\ell(p)-1} \{K_n(x_{j+1,p})-K_n(x_{j,p})\} g(u,c_{j,p}) \rightarrow \Psi_{n,a}(u),
$$

\Bin as $\max\{x_{j+1,p}-x_{j,p}, \ 1\leq j\leq \ell(p)-1\} \rightarrow 0$ as $p\rightarrow +\infty$. We may choose all the $c_{j,p}$ not null from the interior of 
$(x_{j,p}, x_{j+1,p})$ ($x_{j,p}<x_{j+1,p}$). We have

\begin{eqnarray*}
S_p(u)&=&\sum_{j=0}^{\ell(p)-1} \frac{\lambda_{K_n}(]x_{j,p}, x_{j+1,p}])}{c_{j,p}^2} \biggr(e^{ic_{j,p}u}-1-ic_{j,p}u\biggr)\\
&=&\sum_{j=0}^{\ell(p)-1} -i \frac{\lambda_{K_n}(]x_{j,p}, x_{j+1,p}])}{c_{j,p}^2}c_{j,p}u + \frac{\lambda_{K_n}(]x_{j,p}, x_{j+1,p}])}{c_{j,p}^2}
\biggr(e^{ic_{j,p}u}-1\biggr)\\
&=:&\sum_{j=0}^{\ell(p)-1} -i \mu_{j,p} u + \lambda_{j,p} \biggr(e^{ic_{j,p}u}-1\biggr),
\end{eqnarray*}

\Bin with

$$
\lambda_{j,p}=\frac{\lambda_{K_n}(]x_{j,p}, x_{j+1,p}])}{c_{j,p}^2} \ and \ \mu_{j,p}=\frac{\lambda_{K_n}(]x_{j,p}, x_{j+1,p}])}{c_{j,p}^2}c_{j,p}.
$$

\Bin We clearly see that $\exp(S_p)$ is the product of Poisson type \textit{cha.f} converging to $\exp(\Psi_n)$ as $p\rightarrow +\infty$ and $a\rightarrow +\infty$. But we also have that $\exp(\Psi_n)$ is continuous. So by the L\'evy continuity theorem, $\exp(\Psi_n)$ is a \textit{cha.f} and it is \textit{idecomp} by 
Property 4 (see page \pageref{pageprop4}).\\

\Ni  Let us study the differentiability of $\Psi_n$. We have

$$
\left|\frac{\partial g(u,x)}{\partial x}\right|=\left|\frac{ix(e^{iux}-1}{x^2})\right|\leq \frac{2}{|x|} \in \mathcal{L}^1({K_n}),
$$

\Bin and hence

$$
\Psi_n^{\prime}(u)=\int \frac{ix(e^{iux}-1}{x^2}) \ dK_n(x) \ and \ \Psi_n^{\prime}(0)=0.
$$

\Bin Also

$$
\left| \frac{\partial^2 g(u,x)}{\partial^2 x}\right|=\left|\frac{-x^2e^{iux}}{x^2}\right|=1  \in \mathcal{L}^1(K_n),
$$

\Bin and hence

$$
\Psi_n^{\prime\prime}(u)=- \int e^{iux} \ dK_n(x) \ and \ \Psi_n^{\prime\prime}(0)=-s_n^2.
$$

\Bin Finally, let $Z_n$ be a \textit{r.v.} with \textit{cha.f} $\exp(\Psi_n)$. The first and second derivatives  of $\exp(\Psi_n)$ are

$$
\Psi_n^{\prime}(u) \exp(\Psi_n(u)) \ and \ \{\Psi_n^{\prime\prime}(u) \exp(\Psi_n(u)) + \left(\Psi_n^{\prime}(u)\right)^2 \exp(\Psi_n(u))\}  
$$

\Bin taking the values

$$
\Psi_n^{\prime}(0) \exp(\Psi_n(0))=0 \ and \ \{\Psi_n^{\prime\prime}(0) \exp(\Psi_n(0)) + \left(\Psi_n^{\prime}(0)\right)^2 \exp(\Psi_n(0))\}=-s_n^2.  
$$

\Bin We conclude that $\mathbb{E}Z_n=0$ and $\mathbb{V}ar(Z_n)=s_n^2$. The relation 

$$
\Psi_n^{\prime\prime}(u)=- \int e^{iux} \ dK_n(x) \ \ \ (C)
$$

\Bin shows that $\Psi_n^{\prime\prime}(u)$ characterizes $K_n$ and vice-versa. Now, for two functions $\Psi_n$ and $\Phi_n$ such that $\Psi_n^{\prime \prime}=\Phi_n^{\prime \prime}$  with 
$\Psi_n(0)=\Phi_n(0)=0$ and $\Psi_n^{\prime}(0)=\Phi_n^{\prime}(0)=0$,  we have

$$
\forall u \in \mathbb{R}, \ \ \Psi_n^{\prime}(u)=\Phi_n^{\prime}(u) + d_1,
$$

\Bin and by applying this for $u=0$, we get $d_1=0$. Next, we have

$$
\forall u \in \mathbb{R}, \ \ \Psi_n(u)=\Phi_n(u) + d_2,
$$

\Bin and by applying this for $u=0$, we get $d_2=0$. So $\Psi=\Phi$ and  we have the following fact.

\begin{fact} \ \label{charact_01}
$K_n$ characterizes $\Psi_n$ and vice-versa. $\blacksquare$\\
\end{fact}



\section{The weak convergence theorem of summands under the \textit{BVH} and the \textit{UAN} Condition} \label{03_stable_GEN_03}

\subsection{The Central limit theorem for centered, independent and square integrable random variables} \label{clt_subsec_01}$ $\\

\Ni We are going to conclude the discussion above to find solutions of the \textit{CLT} problem under the \textit{BVH} and the \textit{UAN} Condition. We will have two studies from which of them we draw a final conclusion.\\

\Ni \textbf{Study (A)}. From Lemma \ref{03_03_lem_02} and from the notations above, we have

\begin{equation*}
\forall \ t \in \mathbb{R}, \ \Psi_{S_n}(t)-\exp(\Psi_n(t)) \rightarrow 0 \ as \ n\rightarrow +\infty.
\end{equation*}

\Bin But $\exp(\Psi_n(\circ))$ is an \textit{idecomp} \textit{cha.f} for any $n\geq 1$ and is linked to

\begin{equation}
\Psi_n(u)=\int \ g(u,x) \ d\lambda_{K_n}(x) \ with \ g(u,x)=\frac{e^{iux}-1-iux}{x^2}, \ x \in \mathbb{R}, \label{HB}
\end{equation}

\Bin where $\lambda_{K_n}$ is the Lebesgue-Stieltjes measure associated with the \textit{df} $K_n$. Now, we are using the weak convergence theory of bounded measures on $\mathbb{R}$ as exposed Chapter 6 in \cite{ips-wcia-ang}.\\

\Ni \textbf{Direct part}. Let us suppose that $\lambda_{K_n}$ pre-weakly converges to some \textit{df} $\lambda_{K}$, i.e., 
(for $C(K)$ standing for the set continuity points of $K$),

$$
\forall x \in C(K), \ K_n(x)  \rightarrow K(x) \ as \ n\rightarrow +\infty.
$$

\Bin By Part (i) of Proposition 37 in Chapter 6 in \cite{ips-wcia-ang}, we have

$$
\lambda_K(\mathbb{R}) \leq \liminf_{n\rightarrow +\infty} \lambda_{K_n}(\mathbb{R})\leq c,
$$

\Bin since, for any $n\geq 1$,

$$
\lambda_{K_n}(\mathbb{R})=\sum_{k=1}^{k(n)} \int x^2 dF_{k,n}(x)=\sum_{k=1}^{k(n)} \mathbb{V}ar(X_{k,n})\leq c
$$

\Bin from the \textit{BVH}. Hence the pre-weak limit $\lambda_K$ is a bounded measure. Now we apply the integral Helly-Bray theorem as in Theorem 30 in Chapter 6 in \cite{ips-wcia-ang} to \eqref{HB} (See above). By \eqref{bound_01}, for any fixed real number $u$, the function $g(u,x)$ (in $s$) in \eqref{HB} is continuous and satisfies $g(\pm \infty)=0$. So by the cited Helly-Bray integral theorem, we have

$$
\forall u\in \mathbb{R}, \ \Psi_n(u) \rightarrow \Psi_K(u)=\int \ g(u,x) \ d\lambda_{K}(x)=:\int \frac{e^{ix}-1-iux}{x^2} \ d\lambda_{K}(x).
$$

\Bin Now, from the expression of $\Psi_K(u)$ and from \eqref{bound_01}, we see that $\Psi_K(u)$ is a parametrized (in $u$) integral and by the dominated convergence theorem,  $\Psi_K(u)$ is continuous and $\Psi_K(0)=0$. Therefore,

$$
\forall u\in \mathbb{R}, \ \exp(\Psi_n(u)) \rightarrow \exp(\Psi_K(u))=:f_{K}(u).
$$

\Bin Since  $f_K(\circ)$ is continuous at $zero$ and $f_K(0)=1$, we get by the L\'evy continuity theorem (See Theorem 11 in Chapter 3 in \cite{ips-wcia-ang}), we conclude that $f_K$ is \textit{cha.f} and by designating by $\mathcal{K}_K$ the probability law associated to the \textit{cha.f} $f_K$, we have

$$
S_n \rightsquigarrow \mathcal{K}_K.
$$

\Bin By \textit{Property 4} (see page \pageref{property_05} above), $\mathcal{K}$ is an \textit{idecomp} probability law, following the fact that each $\exp(\Psi_n(\circ))$, $n\geq 1$, is an \textit{idecomp} \textit{cha.f}. \\

\Ni \textbf{Indirect Part}. Suppose that for some \textit{df} $K_0$,

$$
S_n \rightsquigarrow \mathcal{K}_{K_0},
$$

\Bin where $\mathcal{K}_{K_0}$ is the probability law associated to $K_0$. We are going to use a Prohorov's type argument. By the asymptotic 
tightness theorem (See Theorem 29 in \cite{ips-wcia-ang}), any sub-sequence $(\lambda_{K_{n_j}})_{j\geq 1}$ of $(\lambda_{K_{n}})_{n\geq 1}$ contains a sub-sequence
$(\lambda_{K_{n_{j_{\ell}}}})_{\ell\geq 1}$ pre-weakly converging to some $\lambda_{K^\ast}$. By the direct part,

$$
S_{n_{j_\ell}} \rightsquigarrow \mathcal{K}_{K^\ast},
$$

\Bin where $\mathcal{K}_{K^\ast}$ is associated to a \textit{cha.f} $f_{K^{\ast}}=\exp(\Psi_{K^\ast})$, with

$$
\forall u\in \mathbb{R}, \ \exp(\Psi_{n_{j_\ell}}(u)) \rightarrow \exp(\Psi_{K^\ast}(u))=:f_{K}(u)
$$

\Bin and

$$
\forall u \in \mathbb{R}, \ \Psi_{K^\ast}(u)=\int \frac{e^{iux}-1-iux}{x^2} \ d\lambda_{K^\ast}(x).
$$

\Bin By uniqueness of the weak limit, $\mathcal{K}_{K^\ast}=_d\mathcal{K}_{K_0}$. Then each sub-sequence of $(K_n)_{n\geq 1}$ contains a sub-sequence converging to $K_0$. We conclude that by Prohorov theorem 

$$
K_n \ \rightsquigarrow_{pre} K_0.
$$

\Bin In both parts, $\mathbb{V}ar(K)<+\infty$ and by Fact \ref{charact_01} applied to $K$, we may conclude that $K$ and $\Psi_K$ characterizes one the other.\\

\Ni We conclude as follows.

\begin{theorem} \label{CLTF_01} Under  the \textit{BVH} and the \textit{UAN} Condition for summands of independent, centered and square integrable real valued random variables, we have:\\

\Ni (a) If 

$$
S_n  \rightsquigarrow \mathcal{K},
$$ 

\Bin where $\mathcal{K}$ is a probability law, then is $\mathcal{K}$ is \textit{idecomp}.\\

\Ni (b) For any  \textit{idecomp} probability law $\mathcal{K}$ of a centered and square integrable  random variable $Z$, for which for any $n\geq 1$, there exists
$X_{1,n}$, $\cdots$, $X_{n,n}$ independent and of same law (they are necessarily centered and square integrable) such that

$$
Z=X_{1,n}+\cdots+X_{n,n}=:S_n.
$$

\Bin Then clearly, $\mathcal{K}$ is a weak limit of summands of independent, centered and square integrable real valued random variables under the \textit{BVH} and the \textit{UAN} Condition.\\

\Ni (3) We have, under the \textit{BVH} and the \textit{UAN} Condition,

$$
S_n  \rightsquigarrow \mathcal{K}_K,
$$ 

\Bin for some \textit{df} $K$, if and only if (using the notation stated above)

$$
K_n \rightsquigarrow_{pre} K.
$$

\Bin Moreover 

$$
\Psi_K(u)=\int \frac{e^{iux}-1-iux}{x^2} \ dK
$$

\Bin and $K$ characterize each other, and $\exp(\Psi_K(\circ))$ is the characteristic function of $\mathcal{K}_K$.
\end{theorem}

\Bin \Ni \textbf{Study (B)}. Here, we suppose that  the \textit{VCH} and the \textit{UAN} Condition hold. We begin by remarking that the Comparison Lemma \ref{03_03_lem_01} holds since formula \eqref{boundX_01} (page \pageref{boundX_01}) holds with the use of the \textit{VCH} in Line (L4).\\

\Ni \textbf{Direct part}. Let $K_n \rightsquigarrow K$. In particular $K_n \ \rightsquigarrow_{pre} K$. By the direct part of \textit{Study (A)}, we still have

\begin{equation} \label{w00}
\lambda_K(\mathbb{R}) \leq \liminf_{n\rightarrow +\infty} \lambda_{K_n}(\mathbb{R})=c,
\end{equation}

\noindent and

$$
S_n \rightsquigarrow \mathcal{K}_K.
$$

\Bin Actually, by weak convergence, we exactly have

$$
\lambda_K(\mathbb{R})=\lim_{n\rightarrow +\infty} \lambda_{K_n}(\mathbb{R})=c,
$$

\Bin but this not play any role for the direct part.\\

\Ni \textbf{Indirect part}. Let 

$$
S_n \rightsquigarrow \mathcal{K}_{K},
$$

\Bin for some \textit{df} $K$. By the indirect part of \textit{Study (A)}, we still have

$$
K_n \ \rightarrow_{pre} K.
$$

\Bin Now if, for $Z\sim \mathcal{K}_{K}$ such that $\mathbb{V}ar(Z)=c$, we have that $\lambda_K(\mathbb{R})=\mathbb{V}ar(Z)$ and then

$$
\lim_{n\rightarrow +\infty} \lambda_{K_n}(\mathbb{R})=\lambda_{K}(\mathbb{R}) \ \  and \ \ K_n \ \rightarrow_{pre} K.
$$

\Bin By Proposition 37 in Chapter 6 in \cite{ips-wcia-ang}, we conclude that $K_n \rightsquigarrow K$ \ as $\rightarrow+\infty$. \\

\Ni We conclude as follows.\\

\begin{theorem} \label{CLTF_02} Under the \textit{VCH} and the \textit{UAN} Condition for summands of independent, centered and square integrable real valued random variables, we have the following characterization. If $\mathcal{K}_K$ is associated with a random variable $Z$ such that $\mathbb{V}ar(Z)=c$, where $c$ is the limit in the \textit{VCH}, then we have

$$
S_n  \rightsquigarrow \mathcal{K}_K,
$$ 

\Bin if and only if 

$$
K_n \rightsquigarrow K.
$$
\end{theorem}

\subsection{The Central limit theorem for non-centered, independent and square integrable random variables} \label{clt_subsec_02}$ $\\

\Ni Let us re-conduct all the notations in Subsection \ref{clt_subsec_01}. Let us denote

$$
\biggr(\forall n\geq 1, \ \forall 1\leq k\leq k(n), \ \mathbb{E}X_{k,n}=a_{k,n}\biggr) \ and \ \biggr(\forall n\geq 1, \ \sum_{k=1}^{k(n)} a_{k,n}=a_n\biggr)
$$

\Bin Let us write

$$
\forall n\geq 1, \ \ S_n=(S_n-a_n) + a_n=\sum_{k=1}^{k(n)}(X_{k,n}-a_{k,n}) + a_n =:S_n^{\ast} + a_n.
$$

\Bin Let us  denote by $F_{k,n}^{\ast}$ the \textit{cdf} of $(X_{k,n}-a_{k,n})$ for $n\geq 1$ and $1\leq k \leq k(n)$,

$$
\forall u\in \mathbb{R}, \Psi_{K^\ast}(u)=\int \frac{e^{iux}-1-iux}{x^2} \ dK^{\ast}(x)
$$

\Bin and

$$
\forall n\geq 1, \ \forall u\in \mathbb{R}, \ \ \Psi_{K_n^{\ast}}(u)=\int \frac{e^{iux}-1-iux}{x^2} \ dK_n^{\ast}(x)
$$

\Bin with

$$
\forall n\geq 1, \ \forall x\in \mathbb{R}, \ \ K^{\ast}_n(x)=\int_{-\infty}^{x} y^2  \sum_{k=1}^{k(n)} \ dF_{k,n}^{\ast}(y).
$$

\Bin \textbf{Direct part}. If $K^{\ast}_n(x) \rightsquigarrow_{pre} K^\ast$ and $a_n \rightarrow a$, then

$$
S_n \rightsquigarrow \mathcal{K}_{K^\ast} + a=\mathcal{K}_{0}.
$$ 

\Bin Moreover, the \textit{cha.f} of $\mathcal{K}_{K^\ast}$ is $\exp(\Psi_{K^\ast}(\circ))$ and next the \textit{cha.f} of $\mathcal{K}_{0}$ is

$$
\forall u\in \mathbb{R}, \ \Psi(u)=\exp\biggr(iau + \Psi_{K^\ast}(u)\biggr).
$$

\Bin \textbf{Indirect part}. Suppose that

$$
S_n \rightsquigarrow \mathcal{K}_{0},
$$ 

\Bin where $\mathcal{K}_{0}$ is associated with an \textit{a.s} finite random variable $Z$. Then $b=\limsup_{n\rightarrow+\infty} a_n$ is finite. Otherwise consider a sub-sequence
$a_{n_\ell}\rightarrow +\infty$ as $\ell\rightarrow +\infty$. So $S_{n_\ell}=S_{n_\ell}^{\ast} + a_{n_\ell}$ necessarily weakly converges to $Z$, where by Theorem \ref{CLTF_01}, $S_{n_\ell}^{\ast} \rightsquigarrow Z^{\ast}$, of law $\mathcal{K}_{K^{\ast}}$ and $Z^{\ast}$ finite \textit{a.s}, and hence $Z$ is \textit{a.s} infinite. Hence $b=\limsup_{n\rightarrow+\infty} a_n$ if finite. Now, each sub-sequence of $(a_{n})_{n\geq 1}$ contains a sub-sequence $(a_{n^\prime})_{n^\prime\geq 1}$ is converging to $a$ finite. By the argument given above, $S_{n^\prime}^\ast$ weakly converges to some $\mathcal{K}_{K^\ast}$. By prohorov's criteria,
$S_{n}^{\ast}$ weakly converges to  $\mathcal{K}_{K^\ast}$ and $S_{n}$ weakly converges to $\mathcal{K}_{K^\ast}+a=_d\mathcal{K}_{0}$. The later inequality shows that all converging subsequences of $(a_{n})_{n\geq 1}$ converge to the same number $a$. Finally 

$$
S_n \rightsquigarrow \mathcal{K}_{K^\ast}+a,
$$
 
\Bin with $a_n\rightarrow a$. Let us summarize the discussions as follows.\\

\begin{theorem} \label{CLTF_03} Under the \textit{BVH} and the \textit{UAN} Condition for summands of independent and square integrable real valued random variables, we have the following characterization. Let us denote

$$
\biggr(\forall n\geq 1, \ \forall 1\leq k\leq k(n), \ \mathbb{E}X_{k,n}=a_{k,n}\biggr) \ and \ \biggr(\forall n\geq 1, \ \sum_{k=1}^{k(n)} a_{k,n}=a_n\biggr);
$$

$$
\forall n\geq 1, \ \ S_n=(S_n-a_n) + a_n=\sum_{k=1}^{k(n)}(X_{k,n}-a_{k,n}) + a_n =:S_n^{\ast} + a_n;
$$

\Bin $F_{k,n}^{\ast}$ the \textit{cdf} of $(X_{k,n}-a_{k,n})$ for $n\geq 1$ and $1\leq k \leq k(n)$;

$$
\forall u\in \mathbb{R}, \Psi_{K^\ast}(u)=\int \frac{e^{iux}-1-iux}{x^2} \ dK^{\ast}(x)
$$

\Bin and, finally,

$$
\forall n\geq 1, \ \forall u\in \mathbb{R}, \ \ \Psi_{K^\ast_n}(u)=\int \frac{e^{iux}-1-iux}{x^2} \ dK_n^{\ast}(x)
$$

\Bin with

$$
\forall n\geq 1, \ \forall x\in \mathbb{R}, \ \ K^{\ast}_n(x)=\int_{-\infty}^{x} y^2  \sum_{k=1}^{k(n)} \ dF_n^{\ast}(y).
$$

\Bin We have the following facts.\\

\Ni (i) If $K_n^\ast \rightsquigarrow_{pre} K^{\ast}$ and $a_n\rightarrow a$, then

$$
S_n \rightsquigarrow \mathcal{K}_{K^{\ast}} + a.
$$

\Bin (ii) If

$$
S_n \rightsquigarrow \mathcal{K}_{0}
$$

\Bin where $\mathcal{K}_{0}$ is associated to an \textit{a.s} finite random variable $Z$, then the sequence $(a_n)_{n\geq 1}$ converges to a real number $a$ and

$$
\mathcal{K}_0 =_d \mathcal{K}_{K^{\ast}} + a
$$

\Bin and 

$$
K_n^\ast \rightsquigarrow_{pre} K^{\ast}.
$$

\Bin Moreover if the \textit{VCH} holds at the place of the \textit{BVH} and the variance of $\mathcal{K}_{K^\ast}$ is equal to $c$,  we have

$$
K_n^\ast \rightsquigarrow K^{\ast}
$$

\Bin in both parts (i) and (ii).
\end{theorem}

\section{Characterizations of two important examples} \label{03_Gauss-Poisson}

\Ni The two important limits of Gaussian law and Poisson law are very important. In stochastic analysis, these laws allow to represent some stochastic process into a discontinuous process (Poisson component) and a continuous process (Gaussian part).

\subsection{Gaussian limit} \label{GC}$ $\\

\Ni Let us suppose that the weak limit of the summands $(S_n)_{n\geq 1}$ is the standard Gaussian law

$$
\exp(\Psi_K(u))=\exp(-u^2/2), \ u\in \mathbb{R},
$$

\Bin i.e.

\begin{equation}
\Psi_K(u)=\int \frac{e^{iux}-1-iux}{x^2} \ dK(x)=-\frac{u^2}{2}, \ u\in \mathbb{R}.
\end{equation}

\Bin But, for $\lambda_{K}=\delta_{0}$, that is, $K=1_{[0,+\infty[}$, we have

$$
\int \frac{e^{iux}-1-iux}{x^2} \ d\delta_0(x)=\left[\frac{e^{iux}-1-iux}{x^2}\right]_{x=0}=-\frac{u^2}{2}.
$$

\Bin We are going to rediscover the L\'evy-Lynderberg-Feller (L2F) theorem as stated in \cite{ips-mfpt-ang} (See Theorem 20 in page ...).

\begin{theorem} \label{L2F-Gauss}
Let $S_n=X_{1,n}+\cdots+X_{k(n),n}$ summands of centered and square integrable random variables as denoted above such that
\begin{equation}
\sum_{1\leq k \leq k(n)} \mathbb{V}ar(X_{k,n})= \sum_{1\leq k \leq k(n)} \sigma^2_{k,n}=1. \label{GC0}
\end{equation}

\Ni For $\varepsilon>0$ and $n\geq 1$, let us denote the Lynderberg function as

\begin{equation}
g_n(\varepsilon)=\sum_{k=1}^{k(n)}\int_{(|x|\geq \varepsilon)} x^2 \ dF_{k,n}(x). \label{lyndG}
\end{equation}

\Bin We have the following characterization:

\begin{equation}
S_n \rightsquigarrow \mathcal{N}(0,1) \ as \ n\rightarrow+\infty \label{GC1a}
\end{equation}

\Ni and

\begin{equation}
\max_{1\leq k\leq k(n)} \mathbb{V}ar(X_{k,n})\rightarrow 0 \ as \ n\rightarrow+\infty   \label{GC1b}
\end{equation}

\Bin if and only if, for any $\varepsilon>0$, the following Lynderberg criterion holds: 

\begin{equation}
g_n(\varepsilon)\rightarrow 0 \ as \ n\rightarrow+\infty \label{GC2}.
\end{equation}

\end{theorem}

\Bin \textbf{Proof}. Let us begin by linking the Lynderberg function as  \eqref{lyndG} with \eqref{GC1b}. We have

\begin{eqnarray}
\max_{1\leq k\leq k(n)} \mathbb{V}ar(X_{k,n})&=& \max_{1\leq k\leq k(n)} \int x^2 dF_{k,n}(x) \notag\\
&=& \sum_{1\leq k\leq k(n)} \int x^2 dF_{k,n}(x) \notag\\
&=& \sum_{1\leq k\leq k(n)} \int_{(|x|\leq \varepsilon)} x^2 dF_{k,n}(x) \notag\\
&+& \sum_{1\leq k\leq k(n)} \int_{(|x|> \varepsilon)} x^2 dF_{k,n}(x) \ \ (L3)\notag \\
&=& \varepsilon^2  + g_n(\varepsilon), \label{GC3}
\end{eqnarray}

\Bin where we used \eqref{GC0} in the first summation in Line (L3). By letting $n\rightarrow +\infty$ first and next, by letting 
$\varepsilon \rightarrow 0$, we get that the Lynderberg criterion implies \eqref{GC1b}. We have: \\

\begin{fact} \label{fact_GC1}
The Lynderberg criterion \eqref{GC2} implies \eqref{GC1b}, which in turn implies the \textit{UAN} hypothesis.
\end{fact}

\Ni Now, Let us prove both implications.\\

\Ni \textbf{Direct implication}. Suppose that \eqref{GC1a} and \eqref{GC1b} hold. So the \textit{BVH} (by \eqref{GC0}) and the \textit{UAN} Condition holds
(by Fact \ref{fact_GC1}). Actually the \textit{BVH} \eqref{GC0} is also a \textit{VCH} Conditions. So may apply both Theorems \ref{CLTF_01} and \ref{CLTF_02}. By applying Theorem \ref{CLTF_01}, we have

$$
\forall x\in C(K), \ K_n(x)=\sum_{1\leq k\leq k(n)} \int_{-\infty}^{x} y^2 \ dF_{k,n}(y) \rightarrow 1_{(x\geq 0)} 
$$

\Bin since $1_{(x\geq 0)}$ is the \textit{df} associated with $\delta_{0}$. Any $x>0$ is in $C(K)$ and then

\begin{eqnarray}
&&K_n(x)=\sum_{1\leq k\leq k(n)} \int_{-\infty}^{x} y^2 \ dF_{k,n}(y) \rightarrow 1 \  \notag\\
&& \Leftrightarrow \sum_{1\leq k\leq k(n)} \int y^2 \ dF_{k,n}(y)-\sum_{1\leq k\leq k(n)} \int_{(y>x)} y^2 \ dF_{k,n}(y) \rightarrow 1 \notag\\
&& \Leftrightarrow 1- \sum_{1\leq k\leq k(n)} \int_{(y>x)} y^2 \ dF_{k,n}(y) \rightarrow 1. \notag
\end{eqnarray}

\Bin Hence

$$
\forall x>0, \ g_{n,1}(x):=\sum_{1\leq k\leq k(n)} \int_{(y>x)} y^2 \ dF_{k,n}(y) \rightarrow 0.
$$

\Bin Next, any $x<0$ is in $C(K)$ and then 

\begin{eqnarray} 
K_n(x)&=&\sum_{1\leq k\leq k(n)} \int_{-\infty}^{x} y^2 \ dF_{k,n}(y) \rightarrow 0 \  \notag\\
&\Leftrightarrow& \sum_{1\leq k\leq k(n)} x^2 \lambda_{F_{k,n}}(\{x\}) +  \sum_{1\leq k\leq k(n)} \int_{(y<x)} y^2 \ dF_{k,n}(y)  \rightarrow 0 \notag \\
&\Leftrightarrow& \lambda_{K_n}(\{x\})+ \sum_{1\leq k\leq k(n)} \int_{(y<x)} y^2 \ dF_{k,n}(y)  \rightarrow 0. \label{discountX} 
\end{eqnarray}

\Bin But, by Portmanteau Theorem (see Criterion (vi) of Theorem 2, page 47 in \cite{ips-wcia-ang}), 
$\lambda_{K_{n}}(\{x\})\rightarrow \lambda_{K}(\{x\})$ since $\partial \{x\}=\{x\}$ and hence $\lambda_{K}(\{x\})=K(x)-K(x-0)=0$ since $x \in C(K)$.\\

\Ni Hence

$$
\forall x<0, \ g_{n,2}(x):=\sum_{1\leq k\leq k(n)} \int_{(y<x)} y^2 \ dF_{k,n}(y) \rightarrow 0.
$$

\Bin By putting together the two last results, for any $\varepsilon>0$

\begin{eqnarray*}
g_n(\varepsilon)&=&\sum_{1\leq k\leq k(n)} \int_{(|y|>\varepsilon)} y^2 \ dF_{k,n}(y)\\
&=&\sum_{1\leq k\leq k(n)} \int_{(y>\varepsilon)} y^2 \ dF_{k,n}(y)+\sum_{1\leq k\leq k(n)} \int_{(y<-\varepsilon)} y^2 \ dF_{k,n}(y)\\
&=&g_{n,1}(\varepsilon) + g_{n,2}(-\varepsilon)\\
&\rightarrow& 0 \ as \ n \rightarrow +\infty. \ \ \square
\end{eqnarray*}

\Bin \textbf{Proof of the indirect implication}. Let \eqref{GC2} holds. So, by Fact \ref{fact_GC1}, \eqref{GC1b} holds and then the \textit{UAN} is satisfied and the \textit{BVH} is already satisfied as an hypothesis of the theorem. Still by Theorem \ref{CLTF_01}, \eqref{GC1a} holds whenever

\begin{equation}
\forall x\in C(K), \ K_n(x)=\sum_{1\leq k\leq k(n)} \int_{-\infty}^{x} y^2 \ dF_{k,n}(y) \rightarrow 1_{(x\geq 0)}. \label{GC4}
\end{equation}

\Bin Let us prove \eqref{GC4}, by exploiting \eqref{GC2}. We have $C(F)=(x<0)+(x>0)$. For $x>0$, we have

\begin{eqnarray*}
K_n(x)&=&\sum_{1\leq k\leq k(n)} \int_{(y\leq x)} y^2 \ dF_{k,n}(y) \\
&=&1 - \sum_{1\leq k\leq k(n)} \int_{(y>x)} y^2 \ dF_{k,n}(y) \\
&=&1 - \sum_{1\leq k\leq k(n)} \int_{(|y|>x)} y^2 \ dF_{k,n}(y) \\
&=&1 - g_n(x) \\
&\rightarrow& 1 \ as \ n \rightarrow +\infty. \ \ \square
\end{eqnarray*}

\Bin For $x<0$, we have

\begin{eqnarray}
K_n(x)&=&\sum_{1\leq k\leq k(n)} \int_{(y\leq x)} y^2 \ dF_{k,n}(y) \notag\\
&=&\sum_{1\leq k\leq k(n)} x^2 \lambda_{F_{k,n}}(\{x\})+\sum_{1\leq k\leq k(n)} \int_{(y<x)} y^2 \ dF_{k,n}(y) \notag\\
&=&\sum_{1\leq k\leq k(n)} x^2 \lambda_{F_{k,n}}(\{x\})+\sum_{1\leq k\leq k(n)} \int_{(-y>-x)} y^2 \ dF_{k,n}(y) \notag\\
&=&\sum_{1\leq k\leq k(n)} x^2 \lambda_{F_{k,n}}(\{x\})+\sum_{1\leq k\leq k(n)} \int_{(|y|> -x)} y^2 \ dF_{k,n}(y) \notag\\
&=& \lambda_{K_n}(\{x\})+ \sum_{1\leq k\leq k(n)} \int_{(|y|> -x)} y^2 \ dF_{k,n}(y) \notag\\
&=& \lambda_{K_n}(\{x\}) + g_{n}(-x). \notag
\end{eqnarray}

\Bin Now, by \eqref{GC2}, $g_{n}(-x)\rightarrow 0$ and by using a similar technical in line \eqref{discountX}, $\lambda_{K_n}(\{x\})\rightarrow 0$. \\

\Ni So \eqref{GC4} holds and we have proved \eqref{GC1a} and \eqref{GC1b}. $\blacksquare$

\subsection{Poisson limit} \label{PC}$ $\\

\Ni The searched limit here is a translated Poisson law $\mathcal{P}(b,\lambda)=b+\mathbb{P}(\lambda)$, with $b\in \mathbb{R}$ and $\lambda>0$ of characteristic function

$$
\exp(\Psi_K(u))=\exp(ibu+\lambda (e^{iu}-1))=\exp(i(b+\lambda)u + \lambda (e^{iu}-1-iu)), \ u\in \mathbb{R},
$$ 

\Bin with

$$
\Psi_K(u)=i(b+\lambda)u + \Psi_{K^\ast}(u), \ \Psi_{K^\ast}(u)=\lambda (e^{iu}-1-iu),
$$

\Bin where $\exp(\Psi_{K^\ast}(\circ))$ is the \textit{cha.f} of the centered Poisson law $\mathcal{P}^\ast(\lambda)=(\mathcal{P}(\lambda)-\lambda)$.\\

\Ni Let us state the characterization theorem.\\

\begin{theorem} \label{L2F-Poiss}
Let $S_n=X_{1,n}+\cdots+X_{k(n),n}$ be summands of independent and square integrable random variables. As above, let $a_{k,n}=\mathbb{E}X_{k,n}$ and let 
$F^{\ast}_{k,n}$ be the \textit{cdf} of $X_{k,n}-a_{k,n}$. Let us introduce  Lynderberg-type  functions, for $\varepsilon>0$ and $n\geq 1$, as

\begin{equation}
g_{n,pois}(\varepsilon)=\sum_{k=1}^{k(n)}\int_{|x-1|> \varepsilon} x^2 \ dF^{\ast}_{k,n}(x). \label{lyndP}
\end{equation}

\Ni Suppose that, as $n\rightarrow +\infty$ with $\biggr(\text{\textit{MVP(n)}}=\sum_{1\leq k \leq k(n)} \sigma^2_{k,n}\biggr)$,

\begin{equation}
B_n=\max_{1\leq k \leq k(n)} \sigma^2_{k,n} \rightarrow 0 \ and \ \text{\textit{MVP(n)}} \rightarrow \lambda.  \label{HPOIS}
\end{equation}

\Bin Let $b\in \mathbb{R}$. We have the following characterization.\\

\begin{equation}
S_n \rightsquigarrow \mathcal{P}(b,\lambda) \ as \ n\rightarrow+\infty \label{PC1}
\end{equation}

\Bin if and only if, 

\begin{equation}
\sum_{k=1}^{k(n)} \mathbb{E}(X_{k,n})=a_n \rightarrow a=b+\lambda  \ as \ n\rightarrow+\infty   \label{PC2a}
\end{equation}

\Bin and for any $\varepsilon>0$, the following Lynderberg Poisson-type criterion holds: 

\begin{equation}
 g_{n,pois}(\varepsilon)\rightarrow 0 \ as \ n\rightarrow+\infty. \label{PC2b}
\end{equation}

\end{theorem}

\Bin \textbf{Proof}. Based \eqref{HPOIS}, the \textit{CVH} and the \textit{UAN} condition  hold. We can apply Theorem \ref{CLTF_03}. We study the limit of $\Psi_{K_n^\ast}(u)$, for any $u \in \mathbb{R}$ to 

$$
\Psi_{K^\ast}(u)=\int \frac{e^{iux}-1-iux}{x^2} \ dK^\ast(x)=\lambda (e^{iu}-1-iu).
$$

\Bin Let $\lambda_{K^\ast}=\lambda\delta_{1}$, i.e., $K^\ast(x)=\lambda 1_{(x\geq 1)}$. Thus

$$
\int \frac{e^{iux}-1-iux}{x^2} \ dK^\ast(x)=\lambda\left[\frac{e^{iux}-1-iux}{x^2} \right]_{x=1}=\lambda(e^{iu}-1-iu), \ u\in \mathbb{R}.
$$

\Ni \textbf{Proof of the direct part}. Suppose that \eqref{PC1} holds. Applying Theorem \ref{CLTF_03}, where the probability law limit is associated with an \textit{a.s} finite random variable, leads to 

$$
a_n \rightarrow b+\lambda \ \ and \ \ K_n^\ast \rightsquigarrow K^\ast.
$$

\Bin So \eqref{PC2a} holds. We also have that $K_n \rightsquigarrow K$ means:

$$
\forall x\in C(K^\ast), \ \ K_n^{\ast}(x) \rightarrow \lambda 1_{(x\geq 1)}, \ as \ n\rightarrow +\infty, 
$$

\Bin since $C(K^\ast)=(x<1)+(x>1)$ and $\lambda_{K_n^\ast}(\mathbb{R}) \rightarrow \lambda_{K^\ast}(\mathbb{R})=\lambda$. For $x>1$, we have

\begin{eqnarray*}
&&\sum_{k=1}^{k(n)} \int_{y\leq x} y^2 \ dF_{k,n}^\ast(y) \rightarrow \lambda\\
&&\Leftrightarrow \sum_{k=1}^{k(n)} \int y^2 \ \ dF_{k,n}^\ast(y) + \sum_{k=1}^{k(n)} \int_{(y>x)} y^2 \ dF_{k,n}^\ast(y) \rightarrow \lambda\\
&&\Leftrightarrow \sum_{k=1}^{k(n)} \sigma^2_{k,n} + \sum_{k=1}^{k(n)} \int_{y> x} y^2 \ dF_{k,n}^\ast(y) \rightarrow \lambda,
\end{eqnarray*}

\Bin where we use that $\int y^2 \ dF_{k,n}^\ast(y)=\mathbb{V}ar(X_{k,n}-a_{k,n})=\sigma^2_{k,n}$ in the last line. Hence

\begin{eqnarray*}
&&\sum_{k=1}^{k(n)} \int_{y\leq x} y^2 \ dF_{k,n}^\ast(y) \rightarrow \lambda\\
&&\Leftrightarrow \sum_{k=1}^{k(n)} \int_{y> x} y^2 \ dF_{k,n}^\ast(y) \rightarrow 0. \ \ (L22)
\end{eqnarray*}

\Bin Let $\varepsilon=x-1>0$, (L22) is equivalent to

$$
\sum_{k=1}^{k(n)} \int_{(y-1> x-1)} y^2 \ dF_{k,n}^\ast(y) \rightarrow 0,
$$

\Bin which is

$$
\sum_{k=1}^{k(n)} \int_{(|y-1|> x-1)} y^2 \ dF_{k,n}^\ast(y) \rightarrow 0,
$$

\Bin that is

$$
g_{n,pois}(x-1) \rightarrow 0.
$$

\Bin Next, For $x<1$, we have

\begin{eqnarray*}
&&\sum_{k=1}^{k(n)} \int_{y\leq x} y^2 \ dF_{k,n}^\ast(y) \rightarrow 0\\
&&\Leftrightarrow \sum_{1\leq k\leq k(n)} x^2 \lambda_{F_{k,n}}(\{x\}) + \sum_{k=1}^{k(n)} \int_{(y<x)} y^2 \ dF_{k,n}^\ast(y) \rightarrow 0\\
&&\Leftrightarrow \lambda_{K_n^\ast}(\{x\}) + \sum_{k=1}^{k(n)} \int_{(1-y>1-x)} y^2 \ dF_{k,n}^\ast(y) \rightarrow 0\\
&&\Leftrightarrow \sum_{k=1}^{k(n)} \int_{(|1-y|>1-x)} y^2 \ dF_{k,n}^\ast(y) \rightarrow 0,
\end{eqnarray*}
 
\Bin where we use that $\lambda_{K_n^\ast}(\{x\})\rightarrow 0$ (as shown in line \eqref{discountX} above) in the last line. Hence 

$$
g_{n,pois}(1-x) \rightarrow 0.
$$

\Bin By combining these results, we have for any $\varepsilon>0$, by taking either $x-1=\varepsilon$ (for $x>1$) or $1-x=\varepsilon$ (for $x<1$), we arrive at \eqref{PC2b}.\\

\Ni \textbf{Proof of the indirect implication}. Suppose that \eqref{PC2a} and \eqref{PC2b} are satisfied. Let us exploit \eqref{PC2b}. For $x>1$, 

\begin{eqnarray*}
K_n(x)&=&\sum_{k=1}^{k(n)} \int_{(y\leq x)} y^2 \ dF_{k,n}^\ast(y)\\
&=&\sum_{k=1}^{k(n)} \sigma^2_{k,n}-\sum_{k=1}^{k(n)} \int_{(y-1>x-1)} y^2 \ dF_{k,n}^\ast(y)\\
&=&\sum_{k=1}^{k(n)} \sigma^2_{k,n}- g_{n,pois}(x-1) \rightarrow \lambda.
\end{eqnarray*}

\Bin For $x<1$

\begin{eqnarray*}
K_n(x)&=&\sum_{k=1}^{k(n)} \int_{(y\leq x)} y^2 \ dF_{k,n}^\ast(y)\\
&=&\lambda_{K_n}(\{x\}) + \sum_{k=1}^{k(n)} \int_{(|y-1|>1-x)} y^2 \ dF_{k,n}^\ast(y)\\
&=&\lambda_{K_n}(\{x\}) + g_{n,pois}(1-x). 
\end{eqnarray*}

\Bin So, $\lambda_{K_n}(\{x\})\rightarrow 0$ is shown exactly as in the lines \eqref{discountX} above. Next $g_{n,pois}(1-x)\rightarrow 0$ is Assumption . Hence
$K_n \rightsquigarrow_{pre} \lambda 1_{(\circ \geq 1)}$. To complete the proof, we remark that for any $n\geq 1$, $K_n(-\infty)=0$ and
$K_n(+\infty)=\sum_{k=1}^{k(n)} \sigma^2_{k,n}$. So

$$
\lambda_{K_n}(\mathbb{R})=K_n(+\infty)=\sum_{k=1}^{k(n)} \sigma^2_{k,n} \rightarrow \lambda=\lambda_{K}(\mathbb{R}).
$$

\Bin By Theorem \ref{CLTF_03}, we conclude that \eqref{PC1} holds. $\blacksquare$\\

\section{Conclusion} 

\Ni We hope that we have given a complete exposition of the theory of the weak limits of independent summands of square integrable random variables on $\mathbb{R}$. The next step will be a re-do of the same theory when the existence of the moments, even the first moment, is not required.


\begin{thebibliography}{99}
\bibitem[Billinsgley (1968)]{billingsley} Billingsley, P.(1968). \textit{Convergence of Probability measures}. John Wiley, New-York.




\bibitem[Lo(2017b)]{ips-mestuto-ang} Lo, G. S. (2017) Measure Theory and Integration By and For the Learner.
\text {SPAS Books Series}. Saint-Louis, Senegal - Calgary, Canada.
Doi :  http://dx.doi.org/10.16929/sbs/2016.0005,
ISBN : 978-2-9559183-5-7.\\


\bibitem[Lo \textit{et al.}(2016)]{ips-wcia-ang} Lo, G.S.(2016). Weak Convergence (IA). Sequences of random vectors. \text {SPAS Books Series}.
Saint-Louis, Senegal - Calgary, Canada. Doi : 10.16929/sbs/2016.0001. Arxiv : 1610.05415. ISBN : 978-2-9559183-1-9  


\bibitem[Lo (2018a)]{ips-probelem-ang} Lo, G.S.(2016). A Course on Elementary Probability Theory. SPAS Editions. Saint-Louis, Calgary, Abuja.
Doi : 10.16929/sbs/2016.0003.


\bibitem[Lo (2018b)]{ips-mfpt-ang} Lo, G.S.(2018). \textit{Mathematical Foundation of Probability Theory}. \text {SPAS Books Series}. Saint-Louis, Senegal - Calgary, Canada. Doi : http://dx.doi.org/10.16929/sbs/2016.0008. Arxiv : arxiv.org/pdf/1808.01713






\bibitem[Lo\`eve (1997)]{loeve}  Michel Lo\`{e}ve (1997). \textit{Probability Theory I}. Springer Verlag.  Fourth Edition.\\ 
\end{thebibliography}
\end{document}